\documentstyle[12pt,amssymb,pstricks,fullpage]{amsart}
\newtheorem{Theorem}{Theorem}[section]
\newtheorem{Proposition}[Theorem]{Proposition}

\newtheorem{Lemma}[Theorem]{Lemma}
\newtheorem{Definition}[Theorem]{Definition}
\newtheorem{Example}[Theorem]{Example}

\newtheorem{Fact}[Theorem]{Fact}
\newtheorem{theorem}{Theorem}[section]
\newtheorem{defn}[theorem]{Definition}

\newtheorem{lemma}[theorem]{Lemma}

\newtheorem{fact}[theorem]{Fact}

\newtheorem{eple}[theorem]{Example}
\newtheorem{rmk}[theorem]{Remarks}
\newtheorem{dsc}[theorem]{Discussion}
\newtheorem{nota}[theorem]{Notation}

\newsavebox{\indbin}
\savebox{\indbin}{\begin{picture}(0,0)
\newlength{\gnu}
\settowidth{\gnu}{$\smile$}
\setlength{\unitlength}{.5\gnu}
\put(-1,-.65){$\smile$}
\put(-.25,.1){$|$}
\end{picture}}
\def \indo {\mathop{\smile \hskip -0.9em ^| \ }}

\newcommand{\be}{\begin{enumerate}}
\newcommand{\bd}{\begin{defn}}
\newcommand{\bt}{\begin{theorem}}
\newcommand{\bl}{\begin{lemma}}
\newcommand{\ee}{\end{enumerate}}
\newcommand{\ed}{\end{defn}}
\newcommand{\et}{\end{theorem}}
\newcommand{\el}{\end{lemma}}

\newcommand{\ov}{\overline}
\newcommand{\cd}{\centerdot}

\newcommand{\CP}{{\mathcal P}}

\newcommand{\CL}{{\mathcal L}}

\newcommand{\CM}{{\mathcal M}}

\begin{document}

\title[Constructing the hyperdefinable group]
{Constructing the hyperdefinable group from the group configuration}

\author{Tristram de Piro, Byunghan Kim and Jessica Millar}

\address{Mathematics Department\\ University of Edinburgh\\
James Clerk Maxwell Building\\ Mayfield Road\\ Edinburgh EH9 3JZ\\
UK}
\address{Mathematics Department\\ MIT\\ 77 Massachusetts Avenue\\
Cambridge, MA 02139\\ USA}
\address{MATX 1220\\ Mathematics Department\\ 1984 Mathematics Road\\
 University of British Columbia\\ Vancouver, BC V6T 1Z2\\ Canada}

\email{depiro@@maths.ed.ac.uk} \email{bkim@@math.mit.edu}
\email{jessica@@math.ubc.ca}
\thanks{The authors were supported by the Seggie Brown
research fellowship, an NSF grant and an NSF grant DMS-0102502, in order}
%\thanks{The second author was supported by an NSF grant, and the 
%third author was supported by an NSF grant DMS-0102052}
\begin{abstract}
Under $\CP(4)^-$-amalgamation, we obtain the canonical
hyperdefinable group from the group configuration.
\end{abstract}

\maketitle

The group configuration theorem for stable theories given by Hrushovski
\cite{Hr1}, which  extends Zilber's result for $\omega$-categorical
theories \cite{Z},
 plays a central role in producing deep results in {\em geometric stability
theory} (For a complete exposition of it, see \cite{P}).
For example, it is pivotal in the proof of the dichotomy
theorem for Zariski' structures (See \cite{HrZ}). It is fair to say
 the group configuration theorem is one of  the  foundational theorems in
geometric stability theory and its applications to algebraic
geometry.

The theorem roughly says that one can get the canonical
non-trivial type-definable group from the group configuration, a certain
 geometrical configuration, in stable theories.
The complete generalization of the theorem into the context of simple
theories seemed  unreachable. In their topical paper
\cite{BTW}, Ben-Yaacov, Tomasic and Wagner generalize the group
configuration theorem by obtaining  an invariant group from the
group configuration in simple theories. However the group they
produce does not completely fit into the first-order context.

On the other hand, Kolesnikov in his important
 thesis \cite{Ko1}, categorizes simple
theories by strengthening the type-amalgamation property (the
independence theorem \cite{KP}), along the  lines of
 early suppositions by Shelah \cite{S} and Hrushovski \cite{Hr2}.
These works suggest to us the possibility of  using higher
amalgamation for the group configuration problem. This approach proves
successful, and in this paper we  succeed in getting
the canonical hyperdefinable group from the group configuration
under stronger type-amalgamation in simple theories. The element
of the group is a hyperimaginary, an
 equivalence class of a type-definable equivalence
relation, and the group operation is type-definable,
hence
the group belongs to the domain of
the standard first-order logic.

We assume that the reader is
familiar with basics of simplicity theory \cite{W}.
Throughout the paper, $T$ is a complete simple
theory. We work in a saturated model $\CM$
of $T$ with
hyperimaginaries, and
 $a,b,...$ are  (possibly infinitary)
 hyperimaginaries, $M,N$ are small elementary submodels.
(Note that tuples from $\CM^{eq}$ are also hyperimaginaries).
As usual, $a\equiv_A b\ (a\equiv^L_A b)$ means $a,b$ have the
same type (Lascar strong type, resp.) over $A$. We point
out that usually $bdd(a)$ denotes
 the {\em set} of all {\em countable} hyperimaginaries definable
over $a$ \cite[3.1.7]{W}. Here, depending on the
 context, it can be  either a
 specific
sequence which  linearly orders the set $bdd(a)$; or,  since a
sequence of hyperimaginaries is again a hyperimaginary (of a large
arity), a fixed hyperimaginary interdefinable with the sequence.

We thank Frank O. Wagner for valuable correspondence on our
earlier note which improves  the presentation.

\section{Around the generalized amalgamation property}

The usual amalgamation property (or the independence theorem) for
 Lascar strong types in simple theories is
stated as follows: For $B\indo_A C$ with $A\subseteq B,C$, if $p$
is a Lascar strong type over $A$ and $p_B,p_C$ are nonforking
Lascar strong type extensions of $p$ over $B$,$C$, respectively,
then there is $d\models p_A\cup p_B$ such that $d\indo_A BC$.  We
 call it `3-amalgamation' \cite{Hr4} (rather than 2-amalgamation
\cite{Ko2}) which shall be  compatible with Definition 1.3. Note
that we can think of $B,C$ (after naming $A$) as two vertices of a
base edge of a triangle and $d$ a top vertex, and $p_B=Lstp(d/B),
p_C=Lstp(d/C)$ are the 2 types to be amalgamated. One would expect
higher amalgamation to be a natural generalization of
3-amalgamation, using a tetrahedron and higher dimensional
simplices instead of a triangle.  Indeed, this is the case, but
the following example draws  attention to why we need extra care
in defining the general $n$-amalgamation property.

\begin{Example}
In the random graph $M$ in $\CL=\{R\}$, choose distinct
$a_i,b_i,c_i\in M$
and imaginary
elements $d_i=\{a_i,b_i\}$ ($i=0,1,2$). We can additionally
assume that
$R(a_0,c_0)\wedge R(b_0,c_1)\wedge \neg R(a_0,c_1)\wedge
\neg R(b_0,c_0)$ and $tp(a_0b_0;c_0c_1)=tp(a_1b_1;c_1c_2)
=tp(a_2b_2;c_2c_0)$. Now it follows that $Lstp(d_2/c_0)=Lstp(d_0/c_0),
Lstp(d_0/c_1)=Lstp(d_1/c_1)$ and  $Lstp(d_1/c_2)=Lstp(d_2/c_2)$.
 However it is easy to see  that  $Lstp(d_0/c_0c_1),
Lstp(d_1/c_1c_2), Lstp(d_2/c_2c_0)$ have no common realization.
\end{Example}

In above example, $\{c_0,c_1,c_2\}$ can be considered as
a base triangle, and  $Lstp(d_0/c_0c_1)$,
$Lstp(d_1/c_1c_2)$, $Lstp(d_2/c_2c_0)$
form other 3 triangles  attached
to the base triangle. The example shows that even if the edges of the 3
triangles are compatible over the base vertices,
there is no common vertex joining the 3 triangles.  On the other hand,
due to the nature of the random graph if we only work in
 the home-sort, then any  desired 3 types attached on a base
 triangle with compatible edges will be  realized.
As we want the notion of higher amalgamation to be preserved in
interpreted theories, Kolesinkov suggests, in his revised works
\cite{Ko1} \cite{Ko2}, the following as higher amalgamation which
we call here {\em $K(n)$-amalgamation}. We briefly explain the
notation. In this paper,
 {\em strong type}  indeed means {\em Lascar strong type}. Likewise,
$p\in S_L(A)$ means $p$ is a Lascar strong type over $A$, and for
$B\subseteq A$, $p\lceil_L B$ (or simply $p\lceil B$) denotes
 $Lstp(a/B)$ for any
(some) $a\models p$. Note that for $ q\in S_L(B)$, $q\subseteq p$
means $p\lceil B=q$ or equivalently $p\vdash q$.

\begin{Definition}
\begin{itemize}
\item
We say strong  types $p_i\in S_L(A_i)$ are {\em   compatible over}
$A(\subseteq A_i)$ if each $p_i$ does not fork over $A$ and for
$i,j$, $p_i\lceil_L A_i\cap A_j= p_j\lceil_L A_i\cap A_j$. (Hence
$p_i\lceil_L A= p_j\lceil_L A$). We say  these $A$-compatible
strong types $p_i$  are {\em (generically) amalgamated}
  if  there is $q\in S_L(\bigcup_i A_i)$ nonforking over $A$ such that
$\cup_ip_i\subseteq q$
(i.e. \  $q\vdash p_i$).
\item
We say $T$ has {\em $K(n)$-amalgamation over $B$} if for $B$-independent
$A=\{ a_1, ..., a_n \}$   and any $B$-compatible $p_i\in
S_L(BA_i)$ where $A_i=A\smallsetminus \{a_i\}$ $(i=1,..,n)$,
whenever  for $a\models p_1\lceil_L B(=p_i\lceil_LB)$
$bdd(aB)\subseteq dcl(aB)$, then $p_1\cup...\cup p_n$ is
generically amalgamated. We say $T$ has {\em $K(n)$-amalgamation}
if it has  {\em $K(n)$-amalgamation} over an arbitrary set.
\end{itemize}
\end{Definition}

The  mend
 is that the realizations of strong types need be boundedly
closed over the parameter set, i.e. in above $bdd(aB)\subseteq
dcl(aB)$. Note that $K(2)$-amalgamation  is equivalent to
3-amalgamation  (usual amalgamation), and due to weak elimination
of imaginaries, it can now be seen that the random graph has
$K(n)$-amalgamation for all $n$. Each stable theory has
$K(n)$-amalgamation as well, by stationarity.

\bigskip

However when we use inductive arguments for example, often we need
to mind not only bounded closures of vertices of amalgamated types
but also those of higher dimensional surfaces as well, since after
naming parameters the surface dimension is increasing. Indeed,
there  exists in the literature another notion of amalgamation,
called $\CP(n)^-$-amalgamation, which was introduced by Hrushovski
\cite{Hr2} prior to  Kolesnikov's work. In the notion, above
concern is already taken care of. {\em Moreover differently from
$K(n)$-amalgamation (or the statement of the independence
theorem), the base simplex is  {\em not} regarded as an
embedded parameter, {\em but}  another type to be amalgamated.} We
think this is conceptually more correct and we shall take it to be
{\em $n$-amalgamation}.

\begin{Definition}
Let $I=\CP(n)^-(=\CP(n)\setminus \{n\})$, ordered by inclusion.
Let $(\{A_i\}_{i\in I},\{\pi^i_j\}_{i\leq j\in I})$ be a directed
family. Namely, each $\pi^i_j:A_i\to A_j$ is an elementary map
between the two sets, and for $i\leq j\leq k\in I$, $\pi^j_k\circ
\pi^i_j=\pi^i_k$ and $\pi^i_i=id$. We say $T$ has  {\em
$\CP(n)^-$-amalgamation}, or simply {\em $n$-amalgamation} if
whenever for any $u\in I$,
 \be
 \item
 $\{\pi^{\{i\}}_u(A_{\{i\}}): i\in u \}$ is
$\pi^{\emptyset}_u(A_{\emptyset})$-independent,
 \item
$A_u=bdd(\cup_{i\in u}\pi^{\{i\}}_u(A_{\{i\}}))$,
 \ee
then  we can extend the direct family to the one indexed by
$\CP(n)$ (by finding $A_n$ and $\pi^j_n$) so that (1),(2) hold for
$n$ too. We say $T$ has  {\em $\CP(n)^-$-amalgamation ($n$-amalgamation) over $A$},
if $A_{\emptyset}=bdd(A)$.
\end{Definition}

Since the definition is not transparent
 to conceptualize with the above
notation, we give a rewritten definition as in \cite{CHr} or
\cite{Hr3}. Recall that when we say a hyperimaginary $b=\bar a/E$
realizes  a type $r$ over $d=\bar c/F$, we mean  $r=r(\bar x)$ is
a (real) type such that i) $r(\bar a)$; ii)  whenever $E(\bar
e,\bar e')$, then $r(\bar e)$ iff $r(\bar e')$; iii) $r(\bar a')$
if $\bar a'/E=f(b)$ for some   $d$-automorphism $f$.
If additionally the converse of iii) holds, we call $r$ a complete
type of $b$ over $d$.

\begin{Definition}
We say $T$ has   {\em $n$-complete amalgamation over a set $B$} if
the following holds: Let $W$ be a collection of subsets of
$\{1,...,n\}=u_n$, closed under subsets. For each $w\in W$,
complete type $r_w(x_w)$ over $B$ is given where $x_w$ is possibly
an infinite set of variables. Suppose that

(1) for $w\subseteq w'$, $x_w\subseteq x_{w'}$ and
 $r_w\subseteq r_{w'}$.

\noindent Moreover for any $a_w\models r_w$,

(2) $\{a_{\{i\}}|i\in w\}$ is $B$-independent,

(3) $a_w$ is as a set $bdd(\cup_{i\in w} a_{\{i\}}B)$ (and the map
$a_w\to x_w$ is a bijection).

\noindent Then there is a complete type $r_{u_n}(x_{u_n})$ over
$B$ such that (1),(2),(3) hold for all $w\in W\cup \{u_n\}$. We
say $T$ has   {\em $n$-complete amalgamation ($n$-CA)} if it has
$n$-complete amalgamation over any set.
\end{Definition}

We leave  the reader to show that  $T$ has    $n$-CA
over $B$ iff $T$ has
$m$-amalgamation over $B$ for all $m\leq n$.
The following can be freely
used: For $B$-independent $A=\{a_1,...,a_n\}$, $\{A_w|w\in
\CP(u_n)\}$ is a partition of $bdd(BA)$, where $A_w=bdd(\cup
a_iB|i\in w) \setminus \bigcup_ {v\in \CP(w)^-}bdd(\cup a_iB|i\in
v)$.  For example $n=2$, $\{bdd(B), bdd(a_1B)\setminus bdd(B),
bdd(a_2B)\setminus bdd(B), bdd(Ba_1a_2)\setminus (bdd(a_1B)\cup
bdd(a_2B))\}$ is a partition of $bdd(Ba_1a_2)$, since using the fact that
$a_1\indo_B a_2$, we have $bdd(a_1B)\cap bdd(a_2B)=bdd(B)$. It also
follows in  1.4, for $v,w\in W$, $x_{v}\cap x_{w}= x_{v\cap w}$.

\medskip

Any simple $T$ has $\CP(3)^-$-amalgamation due to usual
amalgamation, and we
 shall see that 4-amalgamation implies $K(3)$-amalgamation (1.8).
For each $n>2$, there is a simple theory having $n$-CA but not having
$(n+1)$-CA over {\em any set} \cite{Ko2}. (The example also shows
 $n$-amalgamation does not necessarily imply $k$-amalgamaion for $k<n$.) 
All stable theories have $n$-amalgamation over a model (1.6).
 Many important simple structures also have $n$-CA for all $n$
such as the random
graph (1.6), every PAC-structure  (over some parameter)
\cite{Hr3}, and ACFA \cite{CHr}.

In the recent work \cite{KKoT}, corrections of terminologies in
\cite{Ko1}\cite{Ko2} in regard to $n$-CA are made.  For instance,
the definition of {\em $K(n)$-simplicity} is presented in terms of
an {\em infinite} Morley sequence. Kolesnikov's ideas in
\cite{Ko1} go through to show the equivalence of $K(2)$-simplicity
and 4-amalgamation. (The equivalence of $K(1)$-simplicity and
3-amalgamation is the way of proving {\em the independence
theorem} \cite{KP}.) Hence it is naturally conjectured that $T$
being $K(n)$-simple and $T$ having
 $(n+2)$-CA are equivalent, for $n>2$. However surprisingly,
 counterexamples are constructed. Then, the revised concept
 of {\em $n$-simplicity} (implying $K(n)$-simplicity)
  defined via a {\em finite} Morley
 sequence is shown to be equivalent to $(n+2)$-CA for every $n$.

\medskip

The lemma 1.5 and 1.6.1,2  below essentially come from the proof
of the generalized independence theorem  \cite{CHr}. We thank Zoe
Chatzidakis for her explanation.

\begin{Lemma}
Let $T$ be stable.

\be
\item
Suppose that for a set $C$, whenever $a\indo_C b$ with
$b=b_1\cup...\cup b_n$, then $dcl(acl(ab_1C)...acl(ab_nC))\cap
acl(bC)= dcl(acl(b_1C)...acl(b_nC))$\ \ ($\sharp$). Then the
following are satisfied. \be
\item
$tp(acl(ab_1C)...acl(ab_nC)/acl(b_1C)...acl(b_nC))$ is stationary.
\item
Let $A=\{a_1,...a_n\}$, $B=\{c_1,...,c_n\}$ be $C$-independent,
respectively. For $1\leq i\leq n$, let $v_i=\{1,...,n\}\setminus \{i\}$.
Now given $k\leq n$, assume there is a bijective  map
\begin{center}
$h:\cup_{1\leq i\leq k}acl(a_{v_i}C)\to \cup_{1\leq i\leq k}acl(c_{v_i}C)$
\end{center}
 where $a_{v_i}=\{a_j|j\in v_i\}$ such that $h(a_i)=c_i$,
$h\lceil C=id$ and,
 for each $v_i$, $h\lceil acl(a_{v_i}C)$
is elementary. Then $h$ is $C$-elementary. \ee
\item
In fact, the condition ($\sharp$) holds when the set $C$ is a
universe of a model $M$. (Hence (1)(a),(b) also are true over a
model.) \ee
\end{Lemma}

\begin{pf}
We can safely assume  $a,a_i, b_i,c_i$ are finite tuples from
$\CM^{eq}=\CM$.

(1)(a) is immediate from ($\sharp$). (Recall that
$cb(c/d)\subseteq dcl(cd)\cap acl(d)$.)

(1)(b) For $1\leq i\leq k$, let $h_i=h\lceil acl(a_{v_i}C)$. Then let
$h^j=h_1\cup...\cup h_j$ ($h^k=h$) and
$D^j_a=dom(h^j)=\cup_{i=1,...,j}acl(a_{v_i}C)$ and
$D^j_c=ran(h^j)=\cup_{i=1,...,j}acl(c_{v_i}C)$. For induction,
assume $h^{j-1}$ is elementary $(1<j)$. We shall show $h^j$ is
elementary too. Now for each $i<j$ let $w_i= v_j\cap v_i$. Then
$a_{v_i}=\{a_j\}\cup a_{w_i}$.  Now since $h_j$ is elementary,
there  is an automorphism $\hat h_j$ extending $h_j$.  Then by
induction, $\hat h(D^{j-1}_a)$ and $D^{j-1}_c$ have the same type
via $h^{j-1}\circ \hat h^{-1}$, in particular have the same type
over the set $\cup_{i=1,...,j-1}acl(c_{w_i}C)$ fixed by
$h^{j-1}\circ \hat h^{-1}$.  Note now that $c_j\indo_C
c_{w_1}...c_{w_{j-1}}$ and $c_{v_i}=\{c_j\}\cup c_{w_i}$. Hence we
can apply (1)(a) to conclude that $\hat h(D^{j-1}_a)$ and
$D^{j-1}_c=\cup_{i=1,...,j-1}acl(c_{v_i}C)$ also have the same
type over $acl(c_{v_j}C)$, i.e. there is an elementary map $g$
sending $\hat h(D^{j-1}_a)$ to $D^{j-1}_c$ fixing
$acl(c_{v_j}C)=ran(h_j)$. Therefore it follows $h^j(\subseteq
g\circ \hat h)$ is elementary.

(2)\   It suffices to show for $e\in
dcl(acl(ab_1M)...acl(ab_nM))\cap acl(bM)$, $e\in
dcl(acl(b_1M)...acl(b_nM))$. Since  $e\in
dcl(acl(ab_1M)...acl(ab_nM))$, there are $e_1...e_n$ and
 $\CL(M)$-formulas
$\varphi(x;y_1...y_n)$, $\psi_i(y_i,zw_i)$ with
$\varphi(e;e_1...e_n)$, $\psi_i(e_i,ab_i)$ such that  $\models
\varphi(u;v)$ implies  $u$ is definable over $vM$, and
$\psi_i(u',v')$ implies $u'$ is algebraic over $v'M$. Therefore
\begin{center}
$\models \exists y_1...y_n
 (\varphi(e,y_1...y_n)\wedge \bigwedge_i\psi_i(y_i,ab_i))$.
\end{center}
Now since $e\in acl(bM)$, $a\indo_M eb$ and so $tp(a/Meb)$ is a
coheir extension of $tp(a/M)$. Thus we have  $m\in M$ such that
\begin{center}
$\models \exists y_1...y_n
 (\varphi(e,y_1...y_n)\wedge \bigwedge_i\psi_i(y_i,mb_i))$.
\end{center}
Hence $e\in dcl(acl(b_1M)...acl(b_kM))$.
\end{pf}

\begin{Proposition}
\be
\item
Let $T$ be stable. If a set $C$ satisfies ($\sharp$) in 1.5.1,
then for each $n$, $T$ has  $n$-CA over $C$.
\item
All stable theories have $n$-CA over a model.
\item
The random graph has $n$-CA over any set. \ee
\end{Proposition}

\begin{pf}
(1) In a stable theory $T$ we can work in $\CM^{eq}$ and
substitute algebraic closures for bounded closures. We use the
notation in 1.4. It suffices to show the case $W= \CP(u_n)^-$ with
the corresponding types $r_w(x_w) (w\in W)$. Again for $1\leq
i<k\leq n$, let $v_i=\{1,...,n\}\setminus \{i\}$ and $w_i= v_k\cap
v_i$. We shall show that $\cup_{1\leq i\leq n}r_{v_i}$ is
consistent and realized by $\cup_{1\leq i\leq n}a_{v_i}$ such that
$\{a_{\{1\}},...,a_{\{n\}}\}$ is $B$-independent. (Then the type
of its algebraic closure over $B$ extending $\cup_{1\leq i\leq
n}r_{v_i}$ is the  desired $r_{u_n}(x_{u_n})$.) Now due to usual
amalgamation there is $a_{v_1}a_{v_2}\models r_{v_1}\cup r_{v_2}$
such that $\{a_{\{1\}},...,a_{\{n\}}\}$ is $B$-independent. Then
for induction, assume that there is $a_{v_1}...a_{v_{k-1}}\models
r_{v_1}\cup...\cup r_{v_{k-1}}$ $(2<k)$ such that
$a_{v_1}...a_{v_{k-1}}$ extends $a_{\{1\}},...,a_{\{n\}}$. Now let
$b_{v_k}\models r_{v_k}$. Then there is a map $h:\cup_{1\leq i<
k}b_{w_i}\to \cup_{1\leq i< k}a_{w_i}$ such that $h$ sends
$b_{w_i}$ to $a_{w_i}$.
 Hence by 1.5.1(b),
$h$ is elementary and hence  extended to an automorphism $\hat h$.
Then we have $a_{v_{k}}=\hat h(b_{v_{k}})\models r_{v_k}$. Now
then  $a_{v_1}...a_{v_{k}}$  realizes $r_{v_1}\cup...\cup
r_{v_{k}}$ if  for $y=x_{v_k}\cap(x_{v_1}\cup ... \cup
x_{v_{k-1}})$, $a_{v_1}...a_{v_{k-1}}\lceil y=a_{v_k}\lceil y$.
But this clearly holds since from the remark after 1.4,
$y=x_{w_1}\cup ...\cup x_{w_{k-1}}$. This finishes the proof of
(1).

(2) It follows from 1.5.2 and (1) above.

(3) Note that for the random graph $\CM=(\bar M,R)$,
we can work in $\CM^{eq}$ and substitute algebraic closures for
bounded closures. Now since
the random graph  has weak elimination of imaginaries, for any
$A$ there is $A'$ in the home sort $\bar M$
such that $acl(A)=dcl(A')$.
Hence when we check $n$-CA of 1.4, we can assume each $r_w$ is a
type of a set in $\bar M$. Then in $\bar M$, since
$tp(A/B)$ is determined
by equality and $R$ relations of pairs in $A\cup B$,
due to randomness of $R$ we clearly have
the desired  unifying type of a set in $\bar M$.
\end{pf}

However, there  is a stable theory which does not even have
4-amalgamation over an algebraically closed set. We thank Ehud
Hrushovski for supplying us with this  example.

\begin{Example}
Let $A$ be an infinite set with $[A]^2=\{\{a,b\}| a,b\in A, a\ne
b\}$, and let $B=[A]^2\times \{0,1\}$ where $\{0,1\}=
{\mathbb{Z}}/ 2 {\mathbb{Z}}$. Also let $E\subseteq A\times [A]^2$
be a membership relation, and let $P$ be a subset of $B^3$  such
that $((w_1,\delta_1)(w_2,\delta_2)(w_3,\delta_3))\in P$ iff there
are distinct $a_1,a_2,a_3\in A$ such that for
$\{i,j,k\}=\{1,2,3\}$, $w_i=\{a_j,a_k\}$, and
$\delta_1+\delta_2+\delta_3=0$.
 Now let $M$ be a model with the 3-sorted
 universe $A, [A]^2, B$ equipped with relations
$E, P$ and the projection $f:B\to [A]^2$.
Then since $M$ is a reduct of
$(A,{\mathbb{Z}}/ 2 {\mathbb{Z}})^{eq}$, $M$ is stable.
We work in $M$ and show $M$ does not have
$\CP(4)^-$-amalgamation. Note first that
$dcl(\emptyset)=acl(\emptyset)$, and for $a\in A$,
$dcl(a)=acl(a)$.
Now choose distinct $a_1,a_2,a_3,a_4\in A$.
For $\{i,j,k\}\subseteq \{1,2,3,4\}$,
fix an enumeration
$\overline{a_ia_j}=(b_{ij},...)$ of $acl(a_ia_j)$ where
$b_{ij}=(\{a_i,a_j\}, \delta)\in B=[A]^2\times \{0,1\}$. Let
$r_{ij}(x_{ij})=tp(\overline{a_ia_j})$, and let $x^1_{ij}$
be the variable for $b_{ij}$. Note that
 $b_{ij}=(\{a_i,a_j\}, \delta)$ and $b'_{ij}=
(\{a_i,a_j\}, \delta+1)$ have the same type over $a_ia_j$.
Hence there is  $(\overline{a_ia_j})'=(b'_{ij},...)$ also realizing
$r_{ij}(x_{ij})$. Therefore we have  complete types
$r_{ijk}(x_{ijk})$, $r'_{ijk}(x'_{ijk})$  both  extending
$r_{ij}(x_{ij})\cup r_{ik}(x_{ik})\cup r_{jk}(x_{jk})$
realized by some enumerations of $acl(a_ia_ja_k)$ such that,
respectively, $P(x^1_{ij}x^1_{ik}x^1_{jk})\in r_{ijk}$ where as
$\neg P(x^1_{ij}x^1_{ik}x^1_{jk})\in r'_{ijk}$.
Then it is easy to see that
$r_{123}\cup r_{124}\cup r_{134}\cup r'_{234}$
is inconsistent.
\end{Example}

In the example,
 $(\{a_2,a_3\},0)\in dcl(acl(a_1a_2)\cup acl(a_1a_3))$,
since $(\{a_2,a_3\},0)$ is a unique solution to
$P((\{a_1,a_2\},0), (\{a_1,a_3\},0), x)$. But
$(\{a_2,a_3\},0)\notin dcl(acl(a_2)\cup acl(a_3))$, i.e.
1.5.1($\sharp$) does not hold over an algebraically closed set. In
\cite{Hr4}, Hrushovski shows that if a stable $T$ eliminates {\em
generalized finite imaginaries} then $T$  has 4-amalgamation.

\smallskip

\begin{Proposition}
If $T$ has 4-amalgamation over $B$, then it has $K(3)$-amalgamation
over $B$.
\end{Proposition}
\begin{pf}
Assume $T$ has 4-amalgamation. Now suppose that $B$-independent
$A=\{ a_1, a_2, a_3\}$   and $B$-compatible $p_i\in S_L(BA_i)$
where $A_i=A\smallsetminus \{a_i\}$ $(i=1,2,3)$ are given. Also
for $d_i\models p_i\lceil_L B$, $bdd(d_iB)\subseteq dcl(d_iB)$\
(*). Now let $r_{\emptyset}(x_{\emptyset})=tp(bdd(B)/B)$ and let
$r_{i}(x_i)=tp(bdd(a_iB)/B)$, $r_4(x_4)=tp(bdd(d_iB)/B)$ extending
$r_{\emptyset}(x_{\emptyset})$. Now we have
$r_{14}(x_{14})=tp(bdd(a_1d_2B)/B)$ extending $r_{1}\cup r_4$
since due to independence $x_1\cap x_0=x_{\emptyset}$. Let
$f_1:bdd(a_1d_2B)\to x_{14}$ be the realization map. Now note that
due to compatibility of types $p_i$, for $i\in {\mathbb{Z}} / 3
{\mathbb{Z}}$, there is an automorphism $h_i$ sending $d_{i+2}$ to
$d_{i+1}$ fixing $bdd(a_iB)$. Then due to (*), $h_{i+2}\circ
h_i\lceil bdd(d_{i+2}B) =h_{i+1}^{-1} \lceil bdd (d_{i+2}B)$ \
(**). Now via $f_1\circ h_1: bdd(a_1d_3B)\to x_{14}$,
$bdd(a_1d_3B)\models r_{14}(x_{14})$. Then there is
$r_{24}(x_{24})=tp(bdd(a_2d_3B)/B)$  extending $r_{2}(x_2)\cup
r_4(x_4)$ since also $x_{14}\cap x_2=x_{\emptyset}$. Thus by the
map $f_2 \circ h_2: bdd(a_2d_1B)\to x_{24}$ where
$f_2:bdd(a_2d_3B)\to x_{24}$, $bdd(a_2d_1B)\models
r_{24}(x_{24})$. Note that $f_2\lceil bdd(d_3B)=f_1\circ h_1\lceil
bdd(d_3B)$.  We too  have $r_{34}(x_{34})=tp(bdd(a_3d_1B)/B)$
extending $r_{3}(x_3)\cup r_4(x_4)$. Let $f_3:bdd(a_3d_1B)\to
x_{34}$. Note again that $f_3\lceil bdd(d_1B)=f_2\circ h_2\lceil
bdd(d_1B)$.  Now  $f_3\circ h_3 :bdd(a_3d_2B)\to x_{34}$ extends
$f_1\lceil bdd(d_2B): bdd(d_2B)\to x_4$ since from  (**), on
$bdd(d_2B)$, $f_3\circ h_3=(f_2\circ h_2)\circ h_3=(f_1\circ
h_1\circ h_2)\circ h_3=f_1$. Therefore $f_1\cup  f_3\circ h_3:
bdd(a_1d_2B)\cup bdd(a_3d_2B)\to r_{14}(x_{14})\cup
r_{34}(x_{34})$ is a well-defined realization  map extending  the
realizations of $r_j(x_j)$ ($j=1,3,4$). Then now it is easy to
find additional types $r_w$ so that they satisfy (1),(2),(3) of
1.4 for $n=4$. Therefore by 4-amalgamation we have $d(\equiv^L_B
d_i)$ such that $\{a_1,a_2,a_3,d\}$ is $B$-independent and the
type of $bdd(a_1a_2a_3dB/B)$ extends types $r_w$. Obviously, $d$
is the generic realization of $p_1\cup p_2\cup p_3$.
\end{pf}

\medskip

The main theme of this paper is finding the canonical group from the
group configuration, a generalization of the group
 configuration theorem of stable theories into the simple context.
We succeed in obtaining the hyperdefinable group from the group
configuration under 4-amalgamation.  What we are going to use is
4-amalgamation over a parameter properly containing a model (See
the proof of 2.6). But as indicated even a stable theory need not
have such a property. Hence to make it work in more general
context, we introduce the notion of {\em model-$n$-CA}, a little
variation of $n$-CA.

\begin{Definition}
We say $T$ has   {\em model-$n$-complete amalgamation} if the
following holds: Let $u_n=\{1,...,n\}$, and
$W_n=\CP(u_{n+1})\setminus \{u_n\}$. Let $W$ be a collection of
subsets of $W_n$, closed under subsets. For each $w\in W$,
complete type $r_w(x_w)$ over a model $M$ is given where $x_w$ is
possibly an infinite set of variables. Suppose that

(1) for $w\subseteq w'$, $x_w\subseteq x_{w'}$ and
 $r_w\subseteq r_{w'}$.

\noindent Moreover for any $a_w\models r_w$,

(2) $\{a_{\{i\}}|i\in w\}$ is $M$-independent,

(3) $a_w$ is as a set $bdd(\cup_{i\in w} a_{\{i\}}M)$ (and the map
$a_w\to x_w$ is a bijection).

\noindent Then there is a complete type $r_{u_{n+1}}(x_{u_{n+1}})$
over $M$ such that (1),(2),(3) hold for all $w\in W\cup
\{u_{n+1}\}$.
\end{Definition}

Each of stability, $(n+1)$-CA over models, or $n$-CA implies model-$n$-CA 
for every $n$.  Model-$n$-CA also holds in aforementioned algebraic 
examples such as  ACFA and PAC-structures. Model-4-CA is the property we
shall use, hence covers the  case that $T$ is stable.

\section{The group configuration}

\begin{Definition}
By {\em a group configuration}  we mean a 6-tuple of
hyperimaginaries $C=(f_1,f_2,f_3,x_1,x_2,x_3)$ over a hyperimaginary
 $e$ such that, for $\{i,j,k\}=\{1,2,3\}$,
\be
\item
$f_i\in bdd(f_j,f_k;e)$,
\item
$x_i\in bdd(f_j,x_k;e)$,
\item
all other triples and all pairs from $C$
are independent over $e$.
\ee
If the group configuration $C=(f_1,f_2,f_3,x_1,x_2,x_3)$ over
$e$ has the property that $bdd(f_i;e)=bdd(cb(x_jx_k/f_ie);e)$, we
call such $C$
{\em a bounded quadrangle}.
If additionally  $x_i,x_j$ are {\em interdefinable}
over $f_ke$, then we call $C$   {\em a definable quadrangle} over $e$.
\end{Definition}

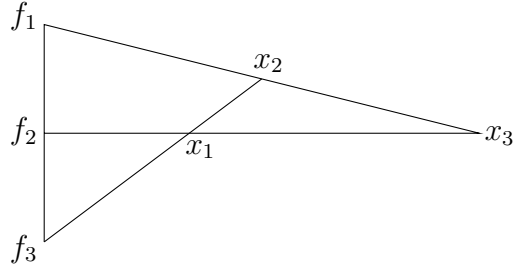
\begin{center}
\psset{unit=0.0012cm}
\begin{pspicture}(5601,2745)(0,-10)
\psset{linewidth=0.4pt}
\psline(375,75)(375,2475)(5175,1275)(375,1275)
\psline(2775,1875)(375,75)
\rput(150,2570){$f_1$}
\rput(150,1300){$f_2$}
\rput(150,0){$f_3$}
\rput(2850,2050){$x_2$}
\rput(5400,1250){$x_3$}
\rput(2100,1100){$x_1$}
\end{pspicture}
\end{center}

\begin{Fact}
\be
\item
If $C=(f_1,f_2,f_3,x_1,x_2,x_3)$ is a group configuration/bounded
quadrangle over $e$ and
$bdd(f_ie)=bdd(f'_ie),bdd(x_ie)=bdd(x'_ie)$, then
$C'=(f'_1,f'_2,f'_3,x'_1,x'_2,x'_3)$ is also a group
configuration/bounded quadrangle over $e$. In this case, we say
$C$ and $C'$ are {\em (boundedly) equivalent} over $e$.
\item
For $C$
a group configuration/bounded quadrangle over $e$ and
$e'\supseteq e$,  if $C\indo_e e'$ then  $C$
also is a group configuration/bounded quadrangle over $e'$.
\item
In a group configuration $(f_1,f_2,f_3,x_1,x_2,x_3)$ over $e$
even if we replace $f_i$ by $f'_i=cb(x_jx_k/ef_i)$ for
$\{i,j,k\}=\{1,2,3\}$,   $(f'_1,f'_2,f'_3,x_1,x_2,x_3)$ is
still a group configuration (hence a bounded quadrangle) over $e$.
\ee
\end{Fact}
\begin{pf} We sketch the proof.
(1) \  Obvious for a group configuration. For a bounded quadrangle
notice that in general $cb(a_1/a_2)$ and
$cb(b_1/b_2)$ are interbounded as far as $a_i,b_i$ are interbounded.\ \
(2) \  Easy.\ \
(3)\  Since $x_ix_j\indo_{f'_k}f_ke$
and  $x_i\indo_{x_jf'_k}f_ke$, $x_i,x_j$ are interbounded over
 $f'_k$\ (*). On the other hand, $x_j\indo_{f_ke} f_i$ implies
$x_ix_j\indo_{f_ke}f_if_j$\ (**), $x_ix_j\indo_{f_kf_ie}f_j$ and
thus $x_ix_jx_k\indo_{f_kf_ie}f_j$ and
$x_ix_jx_k\indo_{f_kf_if_jf'_je}f_j$. Then from (**), it follows
$x_ix_k\indo_{f'_j}f_if_jf_ke$ and from (*)
$x_ix_j\indo_{f'_if'_j}f_if_je$. Therefore $f'_k=cb(x_ix_j/f_ke)
=cb(x_ix_j/f_if_je) \in bdd(f'_if'_j)$. Other independences over
$e$ come easily.
\end{pf}

From now on, assume that a  group configuration over $\hat e=A/\bar E$ is
given. We shall produce the non-trivial canonical hyperdefinable
group from it. By above 2.2.3, we can replace it by  a bounded
quadrangle $C=(\hat f,\hat g, \hat h,\hat a,\hat b,\hat c)$ over a
model $M$ containing $A$. After  naming $M$, we freely assume
that $\emptyset=bdd(\emptyset)$. We further suppose  $\hat f,\hat
g,\hat h,\hat a,\hat b,\hat c$ are all boundedly closed (by
extending each to its bounded closure, if necessary.)
Clearly $C$
still is a bounded quadrangle over $\emptyset$.
 Let $p=tp(\hat f)(=Lstp(\hat f)), q=tp(\hat
g), r=tp(\hat h)$ and let $\Gamma_q(uv)=q(u)\wedge q(v)\wedge
u\indo v$. (Later we shall omit $q$ in $\Gamma_q$.) Now we can
think of $\hat h$ as a multi-valued function such that $dom(\hat
h)= tp(\hat a/\hat h)=Lstp(\hat a/\hat h)$ and $rag(\hat
h)=Lstp(\hat b/\hat h)$. More precisely $b\in k_r(a)$ means
$k_r\models r$, $k_rab\equiv \hat h\hat a\hat b$. Similarly we
write $a\in h_q(c)$, $b\in g_p(c)$ for $h_qca\models \hat g\hat
a\hat c$, $g_pcb\models \hat f \hat b\hat c$, respectively. In the
same way, $b\in dom(f_p)\equiv \exists c (f_pbc\models
 tp(\hat f\hat b\hat c))$, and so on.

\medskip

We say a set $A$ is $n$-independent
if any subset of $A$ having $n$ elements is independent.
Now we define $R=R^{q}$ to be a symmetric type-definable relation
over $\emptyset$
on the set of  independent realizations of $q$ such that
\begin{quote}
$R(fg;f'g')$ iff $f,g,f'g'\models q$, $\{f,g,f',g'\}$
3-independent, and there are $b$ and $a\indo fgf'g'$ such that
$f(a)\cap g(b)\ne\emptyset$, $f'(a)\cap g'(b)\ne\emptyset$.
\end{quote}
It is easy to see that $a\indo fgf'g'$ above can be replaced by
$b\indo fgf'g'$. Similarly, one can define $R^{p}, R^{pq}$ by
replacing $f,g,f'g'\models q$ by $f,g,f'g'\models p$ or
$f,f'\models p,\ g,g'\models q$, respectively.

\begin{Lemma}
\be
\item
If $fg\models \Gamma_q$, and
 $c\in f(a)\cap f(b)$ with $c\indo
fg$, then
\be
\item
$f,g$ are interbounded over
$e:=bc(ba/fg)$, and
\item
 $e\indo f$, $e\indo g$.
\ee
\item
(1) still  holds when we  replace $f,g\models q$ by $f,g\models p$,
or $f\models p, g\models q$.
\item
If $(fg,f'g')\models R$ (or $R^p, R^{pq}$), then  any element in
$\{f,g,f',g'\}$ is in the bounded closure of the other 3 elements.
\ee
\end{Lemma}
\begin{pf}
(1)(a)\ Note that from $ab\indo_e fg$ and that
$a,b$ interbounded over $fg$, it follows that
$a,b$ are interbounded over $e$, too\ (*).
Now from $c\indo_g f$, $ca\indo_g fe$.
Moreover from $c\indo_f g$,
$c\indo_{fe} g$ and then by (*), $ca\indo_{fe}g$. Hence
$g\in bdd(cb(ca/g))\subseteq bdd(fe)$.
By a similar argument, $f\in bdd(ge)$ can be shown too.

(1)(b)\ There are $h_1u_1, h_2u_2$ such that
$h_1fbu_1 c,h_2gau_2 c\models \hat f\hat g\hat a\hat b\hat
c$. Then since $c$ is boundedly closed, by amalgamation
we have
\begin{center}
$hu\models tp(h_1u_1/cbf)\cup tp(h_2u_2/cag) $.
\end{center}
such that $\{u,c,f,g\}$ independent. Then
we have $k,k'$ such that
$hgkauc,hfk'buc\models \hat f\hat g\hat h\hat a\hat b\hat c$.
From $c\indo fgh$, we have $ba\indo_{fg}kk'$. From $b\indo fgh$
and $u,a\in bdd(kk'b)$, it follows
$ba\indo_{kk'}fg$ and thus  $e\in bdd(kk')$\ ($\dag$).
On the other hand,
$f\indo_h g$ implies $k'\indo_h gk$ and  $k'\indo gk$. Hence
$\{g,k,k'\}$ is independent. Similarly $\{f,k,k'\}$ is independent.
Then from ($\dag$), $e\indo f$, $e\indo g$.

(2) Similar to (1).

(3)\ There are $c,c',b$ and $a\indo fgf'g'$ such that
$c\in f(a)\cap g(b)$, $c'\in f'(a)\cap g'(b)$.
Hence $ab\indo_{fg}{f'g'}$ and  $ab\indo_{f'g'}{fg}$.
Therefore $e=cb(ab/fg)$ and $e'=cb(ab/f'g')$ are interbounded\ (**).
From (1)(a),  $f,g$ are interbounded
over $e$, and so are $f',g'$ over $e'$.
Hence it follows from (**),
$g'\in bdd(f'e')= bdd(f'e)\subseteq bdd(f'fg)$
and similarly for the other relations.
\end{pf}

The proof of
2.3.1(b) above is essentially due to   Frank O. Wagner.

\begin{Lemma}
\be
\item
For  $fg,f'g'\models \Gamma_q$, $R(fg,f'g')$ iff there are $b$ and
$a\indo fgf'g'$ such that $f(a)\cap g(b)\ne\emptyset$, $f'(a)\cap
g'(b)\ne\emptyset$
 and $fg\indo_{e}f'g'$ where
$e=cb(ba/fg)$.
\item
Given independent $f,g\models q$, there are $f',g'$ such that
$R(fg,f'g')$.
\item
Above (1)(2) hold for $R^p,R^{pq}$, as well.
 \ee
\end{Lemma}
\begin{pf}
(1)\ ($\Rightarrow$) Note that since $a\indo fgf'g'$,
$ba\indo_{fg}f'g'$, $ba\indo_{f'g'}fg$. Hence
$e,e'=cb(ba/f'g')$ are interbounded. Then due to 2.3.1(a),
$f',g'$ are interbounded over $e$. Now
from $fg\indo f'$, $fg\indo_e f'$, thus $fg\indo_ef'g'$.

 ($\Leftarrow$) Again since $a\indo fgf'g'$,
$e,e'$ are interbounded. Then  from $fg\indo_e f'g'$,
equivalently $fg\indo_{e'} f'g'$, and 2.3.1(b),
$\{f,g, f',g'\}$ is 3-independent.

(2)\  By amalgamation there is $c\in ran(f)\cap ran(g)$ such that
$c\indo fg$. Choose $a\in g^{-1}(c), b\in f^{-1}(c)$.
Then, by the extension axiom,
 we have $f'g'$ such that $\{ab, fg, f'g'\}$ is
$e$-independent where $e=cb(ba/fg)$ and $f'g'\equiv_{abe} fg$.
Then the right hand side of (1) follows easily.

(3)\ Clear.
\end{pf}

The following  lemma is crucial to our argument.

\begin{Lemma}
Let $R(fg,f'g')$.   Namely, $\{f,g,f',g'\}$ 3-independent, and
we can find $d$ and $a\indo fgf'g'$ such that $c\in f(a)\cap g(d)$,
$c'\in f'(a)\cap g'(d)$. Then  there are $h,h'\models p$ and $b$
such that $c\in h(b), c'\in h'(b)$, $b\indo hh'ff'gg'$,
$\{f,h,f',h'\}$ is 3-independent, and $hh'\indo_{ff'}gg'$. It
follows $R^{pq}(hf,h'f')$ and $R^{pq}(hg,h'g')$.
\end{Lemma}

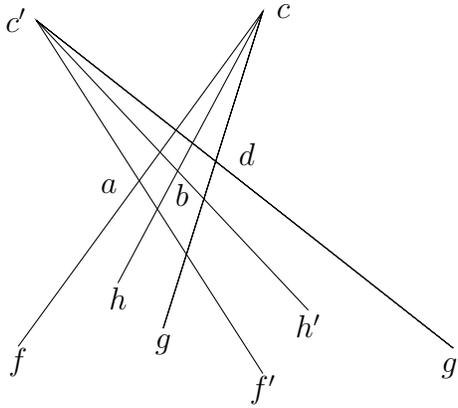
\begin{center}
\psset{unit=0.0012cm}
\begin{pspicture} (5601,4745)(0,-10)
\psset{linewidth=0.4pt}
\psline(375,875)(3075,4575)(1975,1075)(3075,4575)
\psline(3075,575)(575,4475)(5175,855)(575,4475)
\psline(1475,1575)(3075,4575)
\psline(3575,1275)(575,4475)
\rput(375,700){$f$}
\rput(1975,900){$g$}
\rput(1475,1400){$h$}
\rput(3075,400){$f'$}
\rput(5175,680){$g'$}
\rput(3575,1100){$h'$}
\rput(3300,4575){$c$}
\rput(350,4475){$c'$}
\rput(2900,3000){$d$}
\rput(1375,2640){$a$}
\rput(2185,2540){$b$}
\end{pspicture}
\end{center}

\begin{pf}
 Note since
$fac, f'ac'\models \hat g\hat a \hat c$, there are
$k_0,b_0;k_1,b_1$ such that $f'k_0ab_0c', fk_1ab_1c\models \hat
g\hat h\hat a\hat b\hat c$. Then by amalgamation, we have $kb,
hh'$ such that
\begin{center}
$kb\models tp(k_0b_0/af')\cup tp(k_1b_1/af)$, $kb\indo_aff'$ and
$h'f'kabc', hfkabc\models \hat f \hat g\hat h\hat a\hat b\hat c$.
\end{center}
We can further assume that  $hh'\indo_{ff'a}gg'$. Then it follows
from $ab\indo_k ff'$, $b\indo_k ff'hh'$ and $ab\indo_{hh'ff'}gg'$,
$b\indo ff'gg'hh'$. Moreover from $f\indo_kf'$, $hf\indo_kh'f'$,
$\{f,h,f',h'\}$ is 3-independent. Then  since $hh'\indo_{ff'}gg'$,
$\{g,h,g',h'\}$ is 3-independent as well. Therefore $R^{pq}(hf,h'f')$
and  $R^{pq}(hg,h'g')$.
\end{pf}

\medskip

Now define $R'=(R')^q$ by
\begin{quote}
$R'(fg;f'g')$ iff $f,g,f'g'\models q$, $\{f,g,f',g'\}$
3-independent, and for any $a,b$ such that
 $f(a)\cap g(b)\ne\emptyset$ and $a\indo
fg$, there are $a'b'\equiv^L_{fg} ab$  such that $a'\indo
fgf'g'$,   $f'(a')\cap g'(b')\ne\emptyset$.
\end{quote}

Again we also define $(R')^p$, $(R')^{pq}$  by substituting
$f,g,f'g'\models p$ or $f,f\models p, gg'\models q$, respectively.

\medskip

We shall prove that $R$ and $R'$ are equivalent under
4-amalgamation. For the rest of this paper, we assume that $T$ has
4-CA, or more weakly model-4-CA. Note that clearly $R'$ implies
$R$.

\smallskip

\noindent{\bf Notation} For bounded closed sequences $a,b,c$, we
use  $\ov{abc}$ to denote some sequence  of $bdd(abc)$ extending
the orderings of $a,b,c$.

\begin{Theorem}
$R$ and $R'$ are equivalent.
\end{Theorem}
\begin{pf}
It shall show $R$ implies $R'$. Let $R(fg,f'g')$. We  use Lemma
2.5 with the same notation. We then have  $cc',abd,hh'$ such that
$c\in f(a)\cap h(b)\cap g(d)$; $c'\in f'(a)\cap h'(b)\cap g'(d)$;
$b\indo fghf'g'h'$; $f'\in bdd(fhh'), g'\in bdd(hgh')$; and
$\{f,g,h,h'\}$ is independent. Now let $\overline {hf}$,
$\overline {hg}$, $\overline {fg}$, $\overline
{hbc}=bdd(hb)=bdd(hc)$, $\overline {fac}=bdd(fc)=bdd(fa)$,
$\overline {gdc}=bdd(gc)=bdd(gd)$
   be  sequences of boundedly closed sets
 extending the boundedly closed sequences $f,g,h,c,a,b,d$
 (See Notation above 2.6). Since
 $hbc\equiv h'bc'$ we also have
 $\overline{h'bc'}\equiv\overline{hbc}$.

Now, to show   $R(fg,f'g')$, assume there are $a_1,d_1,c_1$ such
that  $c_1\in f(a_1)\cap g(d_1)$ and $c_1\indo fg$. We also have
$\overline {fa_1c_1}\equiv\overline {fac}\equiv \overline
{gd_1c_1}\equiv\overline {gdc}$. Then by 4-amalgamation, there are
$c_2, b_2, a_2, d_2$ such that
\begin{center}
$\overline{fa_2c_2}\ \overline{gd_2c_2}\ \overline{fg}\equiv
\overline{fa_1c_1}\ \overline{gd_1c_1}\ \overline{fg}$;
$\overline{hf}\ \overline{hb_2c_2}\equiv \overline{hf}\
\overline{hbc}$; and $\overline{hg}\ \overline{hb_2c_2}\equiv
\overline{hg}\ \overline{hbc}$\ (*).
\end{center}

 What we are going to amalgamate next are
the following strong types:
\begin{center}
$Lstp(\overline{hb_2c_2}/fg;h)$, $Lstp(\overline{hbc}/h'f;h)$ and
$Lstp(\overline{hbc}/h'g;h)$.
\end{center}
Here the base parameter is $h$ (here is the point where we need
model-4-CA, since indeed the parameter is $Mh$), and each
realization is boundedly closed over the parameter. Each type does
not fork over $h$. Additionally due to (*), it can be seen that
the 3 strong types are $h$-compatible. Hence we have
$\overline{hb_3c_3}$, a generic solution of the types. Moreover we
have $\overline {fhabc}=bdd(f,bc;h)$ extending $\overline {fac}$,
$ \overline {hbc}$.
 Also since $ c'\in bdd(h',bc;h)$
we have $\overline {hh'bcc'}=bdd(h',bc;h)$ extending $\overline
{hbc}$,  $ \overline {h'bc'}$. Similarly  there is $\overline
{hgbdc}=bdd(g,bc;h)$ extending $\overline {hbc}$,  $ \overline
{gdc}$. Note that here  4-amalgamation indeed says that  there exist
elementary maps $\tilde h_1,\tilde h_2,\tilde h_3$ with

\smallskip

$\begin{array}{ll}
\ & dom(\tilde h_1)=bdd(f,bc;h)bdd(h',bc;h)bdd(fh';h),\\
\ & dom(\tilde h_2)=bdd(f,b_2c_2;h)bdd(g,b_2c_2;h)bdd(fg;h),\\
\ & dom(\tilde h_3)=bdd(h',bc;h)bdd(g,bc;h)bdd(h'g;h),
\end{array}$

\smallskip

\noindent fixing $bdd(fh';h)$, $bdd(fg;h)$ $bdd(gh';h)$,
respectively such that  $\tilde h_1(\overline {hbc})=\tilde
h_2(\overline {hb_2c_2})=\tilde h_3(\overline {hbc})=\overline
{hb_3c_3}$ the generic solution. Moreover they are
 compatible with
elementary maps sending $\overline {fhabc}\to \overline
{fha_2b_2c_2}$, $\overline {hgbdc}\to \overline {hgb_2d_2c_2}$ and
$\overline {hh'bcc'} @>{id}>> \overline {hh'bcc'}$. (In
particular, maps $\tilde h_1\lceil bdd(h',bc;h)=\tilde h_3\lceil
bdd(h',bc;h)$.)
 Hence there are $a_3,c'_3,d_3$
such that $a_3=\tilde h_1(a)=\tilde h_2(a_2)$, $c'_3=\tilde
h_1(c')=\tilde h_3(c')$, $d_3=\tilde h_2(d_2)=\tilde h_3(d)$ and
\be
 \item
  $\overline {fha_3b_3c_3}\ \overline {hh'b_3c_3c'_3}
\equiv_{bdd(fh';h)} \overline {fhabc}\ \overline {hh'bcc'}$;
\item
$\overline {fha_3b_3c_3}\ \overline {hgb_3d_3c_3}
\equiv_{bdd(fg;h)} \overline {fha_2b_2c_2}\ \overline
{hgb_2d_2c_2}$;
\item
$\overline {hh'b_3c_3c'_3}\ \overline {hgb_3d_3c_3}
\equiv_{bdd(gh';h)} \overline {hh'bcc'}\ \overline {hgbdc}$.
\ee
Then by (*) and (2),
$a_1d_1\equiv^L_{fg}a_2d_2\equiv^L_{fg} a_3d_3$.
Also since $f'\in bdd(fh'h)$, from (1), $c'_3\in f'(a_3)\cap h'(b_3)$
and $c_3\in f(a_3)$. Note that,  since $R(hg,h'g')$, $g'\in bdd(gh'h)$.
Then, similarly from (3), $c'_3\in g'(d_3)$ and $c_3\in g(d_3)$.
Therefore $R'(fg,f'g')$.
\end{pf}

\section{Type-definability of the transitive closure of $R$}

In this section, we  use $R\equiv R'$ (Theorem
2.6) to prove that,
the transitive closure of $R$ is type-definable. The proof
is similar to the proof in \cite{HKP}  that the transitive
closure of the relation $\sim_1$ forms a hyperimaginary canonical
base, (or the improvement of this proof in \cite[3.3.1]{W}).

\smallskip

Let $\tilde R$ be the transitive closure $R$. We remark that if
both $\{a,b,c,d\}$, $\{a',b',c,d\}$ are 3-independent and
$ab\indo_{cd}a'b'$, then $\{a,b,a',b'\}$ is also 3-independent.

\begin{Lemma}
Suppose that  $R(fg,hk)$, $R(hk,f'g')$ and $fg\indo_{hk} f'g'$. Then
$R(fg,f'g')$ and $fg\indo_{f'g'}hk$.
\end{Lemma}
\begin{pf}
By the previous remark, $\{f,g,f',g'\}$ is 3-independent (*). Now
since $R(fg,hk)$, there are $b$ and $a\indo fghk$
 such that
$f(a)\cap g(b)\ne\emptyset$, $h(a)\cap k(b)\ne\emptyset$
(so $ab\indo_{hk}fg$). Then since
$R'(hk,f'g')$, there are $a'b'\equiv^L_{hk} ab$
 such that $a'\indo
hkf'g'$,   $f'(a')\cap g'(b')\ne\emptyset$.
Hence  by amalgamation, we have
\begin{center}
$a''b''\models Lstp(ab/hk,fg)\cup  Lstp(a'b'/hk,f'g')$,
and $a''b''\indo_{hk}fgf'g'$.
\end{center}
It follows then
$a''\indo fgf'g'$ and $f(a'')\cap g(b'')\ne\emptyset$,
$f'(a'')\cap g'(b'')\ne\emptyset$. This with (*) says
 $R(fg,f'g')$.
It remains to show  $fg\indo_{f'g'}hk$.
Since $g \indo hkf'g'$, $g\indo_{f'g'} hk$. Now by 2.3.3,
$f\in bdd(gf'g')$, and therefore $fg\indo_{f'g'}hk$.
The proof is finished.
\end{pf}

\begin{Theorem}
The following are equivalent.
\be
\item
$\tilde R (\bar f,\bar g)$.
\item
For some $\bar h$, $R(\bar h,\bar f)$ and $R(\bar h,\bar g)$.
\item
For some $\bar h$ with $\bar h\indo_{\bar f}\bar g$ and
$\bar h\indo_{\bar g}\bar f$,
$R(\bar h,\bar f)$ and $R(\bar h,\bar g)$.
\ee
\end{Theorem}
\begin{pf}
It suffices to show (1) implies (3). We prove this by induction on the
length of an $R$-chain. Note that 2.4.2 gives the induction step for
length 0. Now assume that there are $\bar f, \bar f_n, \bar h'$ such that
$\tilde R (\bar f, \bar f_n)$ with the $R$-chain length $n$ and
$R(\bar f_n, \bar h')$.
By the induction hypothesis for $n$, there is $\bar h$ such
$\bar h\indo_{\bar f}\bar f_n$ and
$\bar h\indo_{\bar f_n}\bar f$\ (*),
$R(\bar h,\bar f)$ and $R(\bar h,\bar f_n)$.
 By extension, we can assume
$\bar h\indo_{\bar f\bar f_n}\bar h'$\ (**). Then by (*),
$\bar h\indo_{\bar f_n}\bar h'\bar f$\ (***). In particular,
$\bar h\indo_{\bar f_n}\bar h'$. Hence from the lemma 3.1,
$R(\bar h,\bar h')$ and $\bar h\indo_{\bar h'}\bar f_n$. Then it follows
from (***), $\bar h\indo_{\bar h'}\bar f$.
 Moreover again by (*)(**), we have $\bar h\indo_{\bar f}\bar h'$.
Hence the induction step for $n+1$ is shown.
\end{pf}

\section{The generic group operation on $\Gamma /\tilde R$}

Recall that  in section 2, we define
$\Gamma(xy)=\Gamma_q(xy)=q(x)\wedge q(y)\wedge x\indo y$. Now
since $R$ is symmetric, clearly $\tilde R$ is an equivalence
relation on $\Gamma$. By putting $(\tilde R(\bar x,\bar y)\wedge
\Gamma(\bar x)\wedge \Gamma(\bar y))\vee \bar x=\bar y$, we can
extend $\tilde R$ to  a type-definable equivalence relation on
the whole universe. We shall find the canonical hyperdefinable
group from the hyperdefinable generic  group operation on
 $\Gamma /\tilde R$.
First we state some more
properties of $R$ and $\tilde R$.

\begin{Lemma} Let $\bar f=f_1f_2,\bar g\models \Gamma$, and
let $e=\bar f/\tilde R$.
\be
\item
$\tilde R(\bar f,\bar g)$ and
 $\bar f\indo_e \bar g$
iff
$R(\bar f,\bar g)$.
\item
For $a,b$ such that $f_1(a)\cap f_2(b)\ne \emptyset$ and
$a\indo \bar f$,
$e$ is interbounded with $cb(ab/\bar f)$.
\ee
\end{Lemma}
\begin{pf}
(1)\ ($\Rightarrow$) By 3.2, there is $\bar h$ such that $\bar
h\indo_{\bar f}\bar g$, $\bar h\indo_{\bar g}\bar f$, $R(\bar
h,\bar f)$ and $R(\bar h,\bar g)$. Then since $\bar f\indo_e \bar g$,
${\bar f}\indo_{\bar h}{\bar g}$, and then  by 3.1,\  $R(\bar
f,\bar g)$.

($\Leftarrow$) $R(\bar f,\bar g)$ implies $\tilde R(\bar f,\bar
g)$. By the extension axiom, there is $\bar g'\models tp(\bar
g/e)$ such that $ {\bar g'}\indo_e {\bar g}$ \ (*). Hence $\tilde
R(\bar g,\bar g')$ and   by the proof of ($\Rightarrow$), $R(\bar
g,\bar g')$. Again, by extension, we can assume that
${\bar f}\indo_{e\bar g}\bar g'$.  Hence ${\bar f}\indo_{\bar g}\bar g'$
and then by 3.1, ${\bar f}\indo_{\bar g'}\bar g$.
Therefore by (*),
 ${\bar f}\indo_e{\bar g}$.

(2)\ Let $e_1=cb( ab/\bar f)$. By extension, there is $\bar
f'\models tp(\bar f/abe_1)$ such that $\bar f\indo_{abe_1}\bar
f'$. Hence $\{ab,\bar f,\bar f'\}$ is $e_1$-independent and from
2.4.1, $R(\bar f,\bar f')$ and $e_1=cb( ab/\bar f')$.
 Then by (1), $\bar f\indo_e\bar f'$\ (**).
Let $e_2=cb(\bar f/\bar f')$. Then due to (**), $e_2\in bdd(e)$.
Moreover since $\bar f\indo_{e_2}\bar f'$, $e\indo_{e_2}\bar f'$,
and  $e\in dcl(\bar f')$, $e_2\in bdd(\bar f')$, we have $e\in
bdd(e_2)$. Thus $bdd(e)=bdd(e_2)$. Similarly since $\bar
f\indo_{e_1}\bar f'$, it can be too seen $bdd(e_1)=bdd(e_2)$.
Therefore $bdd(e_1)=bdd(e)$.
\end{pf}

The proof of the following lemma uses 4-CA.

\begin{Lemma}
Let $R(gh,vw)$. Then for any $c\in g(a), d\in v(a)$ with $a\indo
gv$, there are $c'a'd'\equiv^L_{gv}cad$
 and $b'$ such that
$a'\indo ghvw$, $c'\in h(b'),d'\in w(b')$.
\end{Lemma}
\begin{pf}
By 2.5, there are $a_0,b_0,c_0,d_0,t_0$ and $f,u$ such that
 $c_0\in g(a_0)\cap f(t_0)\cap h(b_0)$,
$d_0\in v(a_0)\cap u(t_0)\cap w(b_0)$, $t_0\indo fghuvw$,
$\{f,g,u,v\}$ is 3-independent, and $fu\indo_{gv}hw$. It follows
$R^{pq}(fh,uw)$, $R^{pq}(fg,uv)$ and $gv\indo_{fu}hw$\ (*). Let
$e=bdd(cb(t_0a_0/fg))$. Then from 2.4.3, $e=bdd(cb(t_0a_0/uv))$ and
$fg\indo_{e}uv$. Similarly for
$k=bdd(cb(t_0b_0/fh))=bdd(cb(t_0b_0/uw))$, $fh\indo_{k}uw$. Now
let $\overline {ea_0t_0}=bdd(ea_0t_0)=bdd(ea_0)$, $\overline
{evu}=bdd(ev)$, $\overline {egf}=bdd(eg)$ $\overline
{ga_0c_0}=bdd(ga_0)$, $\overline {va_0d_0}=bdd(va_0)$, $\overline
{gv}$ be sequences of bounded closed sets (See Notation above 2.6).
Note that there are sequences $\overline {gac}$, $\overline {vad}$
such that
\begin{center}
$\overline {ga_0c_0}\equiv \overline {gac}$ and $\overline
{va_0d_0}\equiv\overline {vad}$.
\end{center}
Then by 4-amalgamation, there are $a_1,c_1,d_1,t_1$ such that

\smallskip

$\begin{array}{ll}
(1) & \overline {ga_1c_1}\ \overline {ea_1t_1} \equiv_{ \overline {egf}}
\overline {ga_0c_0}\ \overline {ea_0t_0};\\
(2) & \overline {va_1d_1}\ \overline {ea_1t_1} \equiv_{ \overline {evu}}
\overline {va_0d_0}\ \overline {ea_0t_0};\\
(3) & \overline {ga_1c_1}\ \overline {va_1d_1} \equiv_{ \overline{gv}}
\overline {gac}\ \overline {vad},
\end{array}$

\smallskip

\noindent and $\{a_1,g,v,e\}$ is independent. Then from 2.3.1, it can be seen
so is $\{t_1,f,u,e\}$\ (**). Due to (1)(2), ${ft_0c_0}\equiv
{ft_1c_1}$ and  $ {ut_0d_0}\equiv{ut_1d_1}$. Hence there are
enumerations such that
\begin{center}
$\overline {ft_0c_0}\equiv \overline {ft_1c_1}$ and  $\overline
{ut_0d_0}\equiv\overline {ut_1d_1}$.
\end{center}
Again by 4-CA, we have $c_2,d_2,t_2, b_2$ such that

\smallskip

$\begin{array}{ll}
(4) & \overline {ft_2c_2}\ \overline {kb_2t_2} \equiv_{
\overline {khf}} \overline {ft_0c_0}\ \overline {kb_0t_0};\\
(5) & \overline {ut_2d_2}\ \overline {kb_2t_2} \equiv_{ \overline {kwu}}
\overline {ut_0d_0}\ \overline {kb_0t_0};\\
(6) & \overline {ft_2c_2}\ \overline {ut_2d_2}
\equiv_{ \overline{fu}} \overline {ft_1c_1}\ \overline {ut_1d_1},
\end{array}$

\smallskip

\noindent and $\{t_2,f,u,k\}$ is
independent.  Hence due to  (*),(**),2.3.1(a)
 and (6),
we can apply amalgamation to have
\begin{center}
$d'c't'\models Lstp(d_1c_1t_1/fu;gv)\cup Lstp(d_2c_2t_2/fu;hw)$
\end{center}
such that $t'\indo_{fu}ghvw$\ $(\dag)$. Then there are the desired
$a',b'$ such that

\smallskip

$\begin{array}{ll}
(7) & d'c't'a'\equiv^L_{fguv}d_1c_1t_1a_1,\  \text{and}\ \
d'c't'b'\equiv^L_{fhuw}d_2c_2t_2b_2.
\end{array}$

\smallskip

\noindent Hence it follows from $(\dag)$, $a'\indo ghvw$. Moreover by
(3)(7), $c'a'd'\equiv^L_{gv}cad$; by (4)(7), $c'\in h(b')$; and by
(5)(7), $d'\in w(b')$. The proof is finished.
\end{pf}

\medskip

We are ready to define the promised generic operation on $\Gamma
/\tilde R$.  Let

\smallskip

$\begin{array}{ll}
\cd (x_1y_1,x_2y_2;x_3y_3) := & \exists xyz (\tilde
R(x_1y_1;xy)\wedge \tilde R(x_2y_2;yz)\wedge
  \tilde R(xz;x_3y_3) \wedge \bigwedge_{i=1,2,3} \Gamma(x_iy_i)\\
&\wedge\  x_1y_1/\tilde R\indo x_2y_2/\tilde R\ \wedge
 \{x,y,z\} \text{\ is independent)}.
\end{array}$

\smallskip

\noindent Note that for $x_1y_1,x_2y_2\models \Gamma$,
$x_1y_1/\tilde R\indo x_1y_1/\tilde R$ iff
$\exists x'_1y'_1x'_2y'_2  \tilde R(x_1y_1;x'_1y'_1)
\wedge \tilde R(x_2y_2;x'_2y'_2)\wedge \{ x'_1,y'_1,x'_2,y'_2\}$
 independent. Hence $\cd$ is a partial type over $\emptyset$.

\bigskip

\noindent{\bf Claim 1.} The  relation $\cd$ is
a hyperdefinable partial type over $\emptyset$ such that,
 for any independent
 $e_1=f_1g_1/\tilde R$, $e_2=f_2g_2/\tilde R\in\Gamma  /\tilde R$,
there is $e_3 =fh/\tilde R \in\Gamma  /\tilde R$ such that
$(e_1,e_2;e_3)$ realizes $\cd(x_1y_1,x_2y_2;x_3y_3)$. (See the
explanation above 1.4): It suffices to show  there exist $f,g,h$
such that $e_1=fg / \tilde R$, $e_2=gh / \tilde R$ and $\{f,g,h\}$
is independent.  Now by 2.3.1,  $e_i\indo g_i$,
 $e_i\indo f_i$ and $f_i,g_i$ are interbounded over $e_i$.
Then, by amalgamation, there exists $g\models tp(f_2/e_1)\cup
 tp(g_1/e_2) $ and $\{g,  e_1, e_2\}$ independent. We also have
$f,h$ such that $fg\equiv_{e_1}f_1g_1$, $gh\equiv_{e_1}f_2g_2$.
Then, $\{f,g,h\}$ is independent too.

\medskip

\noindent{\bf Claim 2.}   $e_1\cd e_2=e_3=fh/\tilde R$ does not
depend on the choice of $f,g,h$, i.e. $e_1\cd e_2$ is unique:
Suppose there  are $f',g',h'$ such that $e_1= f'g'/\tilde R $,
$e_2= g'h'/\tilde R$ and $\{f',g',h'\}$ is independent. Then we can
also find independent $\{u, v,w\}$ such that
${u}\indo_{e_1e_2}{fghf'g'h'}$  and $e_1=uv/\tilde R$,
$e_2=vw/\tilde R$. Hence, from 4.1.2,
${uvw}\indo_{e_1e_2}{fghf'g'h'}$, and $\{f,u,e_1,e_2\}$ is
independent\ ($\ddag$).  We shall  prove that $R(fh, uw)$. (Then
the by the same proof, $R(f'h', uw)$ and thus $\tilde
R(fh,f'h')$.) Note  from
 ($\ddag$) and  4.1.1, $R(fg, uv)$  and
$R(gh,vw)$. Hence there are $b$, $a\indo fguv$ and $c\in g(a)\cap
f(b)$, $d\in v(a)\cap u(b)$. Moreover,  by 4.2, we have
$a'b'c'd'$ such that $c'a'd'\equiv^L_{gv}cad$ and $a'\indo ghvw$,
$c'\in g(a')\cap h(b')$, $d'\in v(a')\cap w(b')$. Now again due to
($\ddag$) and 2.3.1, we have $e_1\indo_{gv}e_2$, $fu\indo_{gv}hw$\
$(*)$, and $a\indo_{gv}fu$, $cad\indo_{gv}fu$,
$c'a'd'\indo_{gv}hw$. Hence, by amalgamation, we have
\begin{center}
$c_1a_1d_1\models Lstp(cad/gv;fu)\cup Lstp(c'a'd'/gv;hw)$,
\end{center}
such that $a_1\indo_{gv}fhuw$\ $(**)$. Then we have $b_1,b'_1$
such that
\begin{center}
$c_1a_1d_1b_1\equiv_{fguv}cadb$,
$c_1a_1d_1b'_1\equiv_{ghvw}c'a'd'b'$.
\end{center}
Thus, $c_1\in f(b_1)\cap h(b'_1)$, $d_1\in u(b_1)\cap w(b'_1)$ and
from $(**)$,  $b_1\indo fhuw$.  Moreover from $(*)$ and the remark
above 3.1, $\{f,h,u,w\}$ is 3-independent. Therefore $R(fh,uw)$, as
desired.

\medskip

\noindent{\bf Claim 3.} This generically given group  satisfies
the genericity properties in \cite[4.7.1]{W}: Note that since
$\{f,g,h\}$ independent, it follows for $i=1,2$, $e_i\indo e_1\cd e_2$.
 For  generic associativity, let
$\{k_{1},k_{2},k_{3}\}$ be independent realizations of
$\Gamma/\tilde R$.

\smallskip

\noindent{\bf Subclaim.} There exists independent  $\{h_1,h_2,h_3,h_4\}$
such that
$h_1h_2/\tilde R=k_{1}, h_2h_3/\tilde R=k_{2}$ and $h_3h_4/\tilde R=k_{3}$:
As in the proof of Claim 1,
we can find  ${h_2,h_3,h_4}$ independent such that $h_2h_3/\tilde R=k_{2}$
and $h_3h_4/\tilde R=k_{3}$. Now, for  $h'_1h'_2/\tilde R=k_{1}$,
 amalgamation of $Lstp(h'_2/k_{1})$ and $Lstp(h_2/k_{2}k_{3})$
gives the subclaim.

\smallskip

\noindent Now, by the subclaim, $k_{1}\centerdot k_{2}=h_1h_3/\tilde R$ and
$k_{2}\centerdot k_{3}=h_2h_4/\tilde R$.
Then $k_{1}\centerdot(k_{2}\centerdot k_{3})=h_1h_4/\tilde R$ and
$(k_{1}\centerdot k_{2})\centerdot k_{3}=h_1h_4/\tilde R$ as well.
Finally it can be easily seen that $\cd$ is  generically  surjective.
Hence Claim 3 is verified.

\bigskip

Therefore we have the following;

\begin{Theorem} Given the group configuration,
there exists a canonical hyperdefinable group  and
a definable bijection mapping $\Gamma/\tilde R$  to the generic
types of the group such that $\cd$ is mapped to the group
multiplication generically.
\end{Theorem}

\section{1-based theories}

One application of 4.3 is the following result.
 This extends the theorem \cite[3.23]{DK} that,
in any $1$-based non-trivial
$\omega$-categorical simple $T$, an infinite vector space over some
finite field is
definably recovered in $\CM^{eq}$. Recall that $T$ is {\em non-trivial} if
there are
hyperimaginaries $a_1,a_2,a_3$ and $A$ such that for
$1\leq i<j\leq 3$, $a_i,a_j$ are independent over $A$ whereas
$\{a_1,a_2,a_3\}$ is dependent over $A$.

\smallskip

\begin{theorem}
Suppose that $T$ is $1$-based, non-trivial, having model-4-CA.
Then there is  a hyperdefinable  infinite bounded-by-Abelian group
$V$ over a model $M$ of $SU$-rank 1 generic types. Moreover for
the bounded subgroup $V_0=V\cap bdd(M)$, $V/V_0$ forms a vector
space over the division ring $R$ 
of $bdd(M)$-endomorphisms of $V$ such that 
for $b,a_1,...,a_n\in V$, $b\in bdd(a_1...a_n)$ iff 
$b+V_0=\alpha_1(a_1+V_0)+...+\alpha_n(a_n+V_0)$ for some 
$\alpha_i\in R$.
\end{theorem}

\begin{pf}
By the proof of Lemma 3.22 in \cite{DK}, there exists a non-trivial
 rank-$1$ Lstp $p$ over some model $M$.
For convenience, let $M=\emptyset$ after naming the model.
As $p$ is non-trivial, there
exists $\{a,b,c\}$ realizing $p$ such that ${b,c}$ is independent
and $a\in bdd(b,c)\setminus bdd(b)\cup bdd(c)$. Let $yx$ realize
$tp(ab/c)$
 with $yx\indo_c ab$. Then $dim(ay/bx)=1$ as $y\in bdd(abx)$ and
 $a\indo bx$. Let $z=cb(Lstp(ay/bx))$,
then by $1$-basedness, $z\in bdd(ay)\cap bdd(bx)$.
Moreover, by a straightforward rank calculation, $SU(z)=1$.
This gives a bounded quadrangle $(a,b,c,x,y,z)$.  Now by
 Theorem 4.3, we obtain a hyperdefinable group $G$ over
$\emptyset$ such that the generic types all have $SU$-rank $1$.
The group $G$ is 1-based since the underlying theory is 1-based.
 Now we use the following fact \cite[4.8.4]{W},

\begin{fact}
Suppose $G$ is an $1$-based group hyperdefinable over $\emptyset$
in a simple theory.
Then  for the normal subgroup $G^0_{\emptyset}$,
the smallest $\emptyset$-hyperdefinable subgroup
of bounded index, the commutator subgroup $(G^0_{\emptyset})'$ of
$G^0_{\emptyset}$ has
boundedly many elements and contained in the center of $G^0_{\emptyset}$.
\end{fact}

\noindent Therefore, if we set $G^0=V$, then $V$ is the desired
bounded-by-Abelian hyperdefinable group. Note that by above $V'$
is contained in the normal subgroup $V_0=V\cap bdd(\emptyset)$.
Indeed again from \cite[4.8.18]{W}, the Abelian group $V/V_0$
forms a vector space over a division ring $R$ of
$bdd(\emptyset)$-endomorphisms of $V$, and  dependence in 
$V/V_0$ is given by linear dependence of the vector space.
\end{pf}

\end{document}